\documentclass{ifacconf}

\usepackage{graphicx}      
\usepackage{natbib}        
\usepackage{amsmath,amssymb}
\usepackage{algorithm}
\usepackage{algpseudocode}
\begin{document}
\begin{frontmatter}

\title{Mesh Refinement with Early Termination for Dynamic Feasibility Problems\thanksref{funding}} 

\thanks[funding]{This work was supported by the Science and Technology Facilities Council (STFC) under a doctoral training grant (ST/V506722/1)}

\author[main]{Eduardo M. G. Vila} 
\author[eric]{Eric C. Kerrigan}
\author[paul]{Paul Bruce}

\address[main]{Department of Electrical and Electronic Engineering, Imperial College London, SW7 2AZ, UK (e-mail: emg216@imperial.ac.uk)}
\address[eric]{Department of Electrical and Electronic Engineering, and Department of Aeronautics, Imperial College London, SW7 2AZ, UK (e-mail: e.kerrigan@imperial.ac.uk)}
\address[paul]{Department of Aeronautics, Imperial College London, SW7 2AZ, UK (e-mail: p.bruce@imperial.ac.uk)}

\begin{abstract}                
We propose a novel early-terminating mesh refinement strategy using an integrated residual method to solve dynamic feasibility problems. As a generalization of direct collocation, the integrated residual method is used to approximate an infinite-dimensional problem into a sequence of finite-dimensional optimization subproblems. Each subproblem in the sequence is a finer approximation of the previous. It is shown that these subproblems need not be solved to a high precision; instead, an early termination procedure can determine when mesh refinement should be performed. The new refinement strategy, applied to an inverted pendulum swing-up problem, outperforms a conventional refinement method by up to a factor of three in function evaluations.
\end{abstract}


\end{frontmatter}

\section{Introduction}
Dynamical systems are modeled using differential (and algebraic) equations. The dynamic variables are often constrained by lower and upper bounds, as well as by initial and/or terminal conditions. Such a formulation constitutes a dynamic feasibility problem (DFP) which, in general, may have zero, one, or many feasible solutions. In the latter case, a merit function can be used to find an \emph{optimal} solution, resulting in a dynamic optimization problem (DOP). Many optimal control and system identification problems can be formulated as DOPs. Moreover, a sequence of DFPs may be used to solve a DOP \citep{nie_solving_2022, CDCInteresso}.

In contrast to ordinary optimization problems, where the decision variables are finite-dimensional, DFPs are infinite-dimensional problems: dynamic variables are functions of time (i.e., trajectories). In general, these problems are intractable. Discretization methods are used to \emph{transcribe} the DFP into an approximate, finite-dimensional optimization problem, which is tractable by well-known algorithms.

Within (direct) transcription, the class of \emph{shooting} methods aims to find control variables such that the propagated state variables satisfy the problem constraints. They rely on the sensitivities of a time-marching integrator to iteratively find a solution. Shooting methods often have highly nonlinear relations between control variables and state constraints, resulting in poor convergence rates \citep[Sect. 9.6]{kelly_introduction_2017}. The class of \emph{collocation} methods takes a different approach: both states and controls are discretized and become decision variables of the optimization problem. The dynamic equations are included as constraints, collocated at the discretization points. Direct collocation methods mainly differ in the numerical integration scheme that is used. Examples include trapezoidal and Hermite-Simpson collocation \citep[Sect. 3 and 4]{kelly_introduction_2017}, as well as certain Runge-Kutta and Gauss-Legendre families of collocation methods \citep[Sect. 4.5]{betts_practical_2020}.

Recently, the class of \emph{integrated residual} methods has enjoyed renewed interest \citep{neuenhofen_integral_2020, nie_solving_2022}. Popular among partial-differential equation solvers, integrated residual methods aim to minimize the error on the integration of differential equations. This class of methods can be seen as a generalization of direct collocation, sharing a favorable convergence rate. Belonging to this class is the \emph{least-squares} method, which results in optimization problems of slightly increased dimension; however, the method allows for problems with many algebraic equations and singular arcs to be solved in cases where direct collocation fails \citep{nie_solving_2022}.

In practice, it is more computationally efficient to solve a \emph{sequence} of transcribed optimization problems, rather than a single one. Due to the presence of discontinuities and regions of high gradients, it is favorable to start with a coarse discretization (i.e., fewer mesh points), and progressively refine the mesh towards finer discretizations. The mesh may be refined by increasing the partitioning of the time domain ($h$-methods), and/or by increasing the degree of discretization within each interval ($p$-methods).

State-of-the-art mesh refinement strategies solve the transcribed optimization problems to a high precision. To improve each approximate solution, an estimate of the local discretization error is used to refine the mesh by $h$ and/or $p$ methods. The general principle is to predict the mesh with the least number of points necessary to achieve a desired accuracy. The mesh prediction can use an integer programming technique so as to minimize the maximum local error \citep{betts_mesh_1998}. For problems with smooth solutions, $p$-methods thrive due to their exponential convergence rate \citep{trefethen_spectral_2000}. In the presence of non-smoothnesses, an $h$-method increases convergence by partitioning the time domain into smooth sub-intervals. Thus, a $p$-then-$h$ approach has shown significant computational efficiency \citep{patterson_ph_2015}.

Equally important to \emph{how} the mesh should be refined, is \emph{when} the refinement should take place. For shooting methods, strategies have been proposed where a refinement test is used to terminate each optimization problem early \citep[Sect. 3.3.3]{polak_optimization_1997}. Instead of attempting to predict the ideal mesh prematurely, this approach gradually refines the mesh, as each optimization problem progresses. Refinement tests may be of the \emph{filtering} type, where the early termination criterion is prescribed as a decreasing tolerance. Alternatively, \emph{convergence} type tests terminate each optimization problem once their rates of convergence stagnate, ensuring the convergence rate of the infinite-dimensional problem is achieved \citep{polak_use_1993}.

Our contribution is a progressive mesh refinement strategy for the general class of integrated residual methods. The proposed algorithm simultaneously partitions the mesh ($h$-method), as well as adjusts the integration scheme, ensuring sufficient integration accuracy. We use an early termination test based on the estimated convergence rate which interrupts the optimization problem, and triggers the mesh refinement scheme. The novel use of an early termination test for the class of integrated residual methods results in significant computational savings, when compared with a predictive refinement strategy.

The remainder of this paper is organized as follows. Section~\ref{sec:problem} defines the DFP, and describes the discretization process using an integrated residual method. Section~\ref{sec:refinement} outlines two mesh refinement strategies: progressive refinement (using early termination), and predictive refinement (without early termination). Section \ref{sec:example} applies both strategies to the inverted-pendulum swing-up problem, where the performance is compared. Section \ref{sec:conclusions} suggests future improvements and extensions for the current implementation.

\emph{Notation}: Considering vectors in $\mathbb{R}^n$, $\|\cdot\|$ denotes the Eucledian norm of a vector, and $\langle \cdot, \cdot \rangle$ denotes the inner product of two vectors. Given a set of indices $\mathcal{S} \subseteq \{1, 2, ..., n\}$, $\| \cdot \|_{\mathcal{S}}$ denotes the Euclidean norm of a vector, only over the components in $\mathcal{S}$, and $\langle \cdot, \cdot \rangle_{\mathcal{S}}$ denotes the inner product of two vectors, only over the components in $\mathcal{S}$. $C^0(\Omega)$ denotes the set of continuous functions defined over the domain $\Omega$. Given a multivariate function $f: \mathbb{R}^n \rightarrow \mathbb{R}$, $\nabla f(\cdot)$ denotes the gradient of $f(\cdot)$. The ceiling function $\lceil \cdot \rceil$ yields the smallest integer greater or equal to the argument.

\section{Problem Statement}
\label{sec:problem}
\subsection{Dynamic Feasibility Problem}
Within a continuous-time domain $\mathcal{T} := [t_0, t_f] \subset \mathbb{R}$, we are interested in finding continuous differential dynamic variables $\boldsymbol{x} : \mathcal{T} \rightarrow \mathbb{R}^{n_x}$ (e.g., states of a system) and algebraic dynamic variables $\boldsymbol{u} : \mathcal{T} \rightarrow \mathbb{R}^{n_u}$ (e.g., control inputs) which satisfy a \emph{residuals} function $r : \mathbb{R}^{n_x} \times \mathbb{R}^{n_x} \times \mathbb{R}^{n_u} \rightarrow \mathbb{R}^{n_r}$ (e.g., system dynamics). The dynamic variables are bounded by sets of box constraints $\mathcal{X} \subseteq \mathbb{R}^{n_x}$ and $\mathcal{U} \subseteq \mathbb{R}^{n_u}$, respectively. Moreover, $\boldsymbol{x}(\cdot)$ may be subject to initial conditions $\mathcal{X}_0 \subseteq \mathcal{X}$ as well as terminal conditions $\mathcal{X}_f \subseteq \mathcal{X}$. We formulate the problem as
\begin{equation}
\begin{aligned}
    \text{find} \quad   & \boldsymbol{x}(\cdot), \, \boldsymbol{u}(\cdot),\\
    \text{s.t.} \quad   & r(\dot{\boldsymbol{x}}(t), \boldsymbol{x}(t), \boldsymbol{u}(t)) = 0 \quad &\widetilde{\forall} t \in \mathcal{T}, \\
                        & \boldsymbol{x}(t) \in \mathcal{X}, \, \boldsymbol{u}(t) \in \mathcal{U} \quad &\forall t \in \mathcal{T},\\
                        &\boldsymbol{x}(t_0) \in \mathcal{X}_0, \, \boldsymbol{x}(t_f) \in \mathcal{X}_f,\\
                        &\boldsymbol{x}(\cdot) \in C^0(\mathcal{T}),
\end{aligned}
\tag{\textbf{DFP}}
\label{eq:dfp}
\end{equation}

where $\dot{\boldsymbol{x}} : \mathcal{T} \rightarrow \mathbb{R}^{n_x}$ is the time-derivative of $\boldsymbol{x}(\cdot)$. Because $\boldsymbol{x}(\cdot)$ need only be piece-wise differentiable, the residuals function is enforced `almost everywhere' in the Lebesgue sense ($\widetilde{\forall}t \in \mathcal{T}$). The algebraic variables $\boldsymbol{u}(\cdot)$ are allowed to be discontinous.

We include the dynamic equations implicitly as residuals so as to allow for a broader class of problems that may contain algebraic equations.

\subsection{Integrated Residual Method}
The transcription process begins by reformulating the constraint satisfaction problem \eqref{eq:dfp} into an optimization problem. Using a least-squares approach, the normalized integrated squared residual is
\begin{equation}
    \boldsymbol{f}(\boldsymbol{x}, \boldsymbol{u}) := \frac{1}{n_r(t_f - t_0)}\int_{\mathcal{T}} \| r(\dot{\boldsymbol{x}}(t), \boldsymbol{x}(t), \boldsymbol{u}(t)) \|^2 \, \mathrm{d}t,
\label{eq:irm}
\end{equation}
which should be minimized. To simplify notation, we define the pair $\boldsymbol{z} := (\boldsymbol{x}, \boldsymbol{u})$. The integrated residual $\boldsymbol{f}$ is set as the objective functional of the infinite-dimensional optimization problem
\begin{equation}
     \min_{\boldsymbol{z} \in \boldsymbol{\mathcal{Z}}} \quad \boldsymbol{f}(\boldsymbol{z}),
\tag{\textbf{P}}
\label{eq:p}
\end{equation}
where the constraint set is
\begin{multline}
    \boldsymbol{\mathcal{Z}} := \{ (\boldsymbol{x}, \boldsymbol{u}) \, | \, \boldsymbol{x}(t) \in \mathcal{X}, \, \boldsymbol{u}(t) \in \mathcal{U} \quad \forall t \in \mathcal{T},\\
    \boldsymbol{x}(t_0) \in \mathcal{X}_0, \, \boldsymbol{x}(t_f) \in \mathcal{X}_f, \, \boldsymbol{x}(\cdot) \in C^0(\mathcal{T}) \}.
\end{multline}
Because \eqref{eq:p} is infinite-dimensional, the problem is, in general, intractable. \eqref{eq:p} needs to be discretized by suitable numerical methods so as to achieve an approximate solution.

\subsection{Problem Discretization}
Firstly, the time domain $\mathcal{T}$ is partitioned by $n_{\mathcal{H}}-1$ points $\{ t_1, t_2, ..., t_{n_{\mathcal{H}}-1}\}$, while ensuring
\begin{equation}
    t_0 < t_1 < t_2 < ... < t_{n_{\mathcal{H}}-1} < t_f.
\end{equation}
By change of variables, we define normalized time intervals $\tau_i \in [-1, 1]$ such that
\begin{equation}
    t = \frac{t_i + t_{i-1}}{2} + \frac{t_i - t_{i-1}}{2} \tau_i,
\end{equation} 
where $t \in [t_{i-1}, t_i]$, $\forall i \in \{1, 2, ..., n_{\mathcal{H}}\}$. Additionally, let the normalized length of each interval be
\begin{equation}
    h_i := \frac{t_{i+1} - t_i}{t_f - t_0}.
\end{equation}
The set of interval lengths is $\mathcal{H} := \{h_1, h_2, ..., h_{n_{\mathcal{H}}}\}$. 

\begin{figure}
    \centering
    \includegraphics[width=0.85\columnwidth]{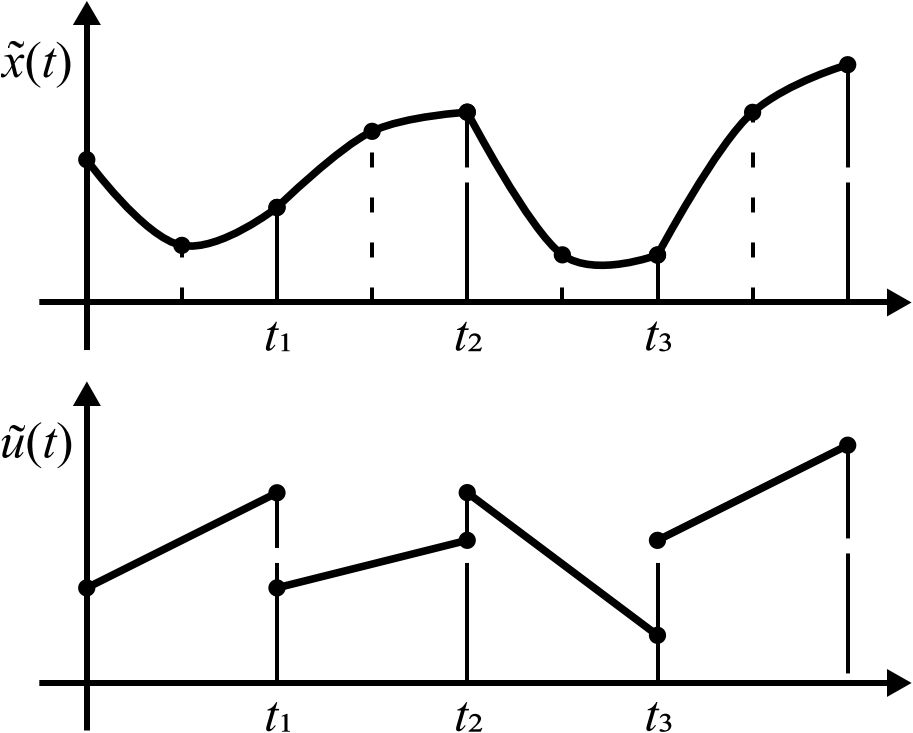}
    \caption{Illustrations of piece-wise Lagrange polynomial approximations of second and first degree, respectively.}
    \label{fig:approximation}
\end{figure}

Secondly, the dynamic variables $\boldsymbol{x}(\cdot)$ and $\boldsymbol{u}(\cdot)$ are approximated by the piece-wise polynomial functions $\tilde{\boldsymbol{x}}(\cdot, \cdot)$ and $\tilde{\boldsymbol{u}}(\cdot, \cdot)$. Within each interval, $\tilde{\boldsymbol{x}}(\cdot, \cdot)$ and $\tilde{\boldsymbol{u}}(\cdot, \cdot)$ are polynomials in $\tau_i$ of degree $n_{\mathcal{P}_x}$ and $n_{\mathcal{P}_u}$, respectively. Between intervals, continuity of the differential variables is enforced by having $\tilde{\boldsymbol{x}}_{i-1}(1, \cdot) = \tilde{\boldsymbol{x}}_i(-1, \cdot)$ $\forall i \in \{2, 3, ..., n_{\mathcal{H}}\}$, as illustrated in Figure \ref{fig:approximation}.

Using Lagrange polynomials in the barycentric form, the approximate dynamic variables $\tilde{\boldsymbol{x}}(\cdot, \cdot)$ and $\tilde{\boldsymbol{u}}(\cdot,\cdot)$ are parameterised by the decision variables $x$ and $u$, respectively. In the $i^{\text{th}}$ sub-interval, the $j^{\text{th}}$ component of $\tilde{\boldsymbol{x}}(\cdot, \cdot)$ is given by
\begin{equation}
    \tilde{\boldsymbol{x}}_j^i(\tau, x) = \frac{\displaystyle \sum_{p=0}^{n_{\mathcal{P}_x}} \frac{w_{p}}{\tau - \tau_{p}} x_{p,j}^i}{\displaystyle \sum_{p=0}^{n_{\mathcal{P}_x}} \frac{w_{p}}{\tau - \tau_{p}}}
\label{eq:bary}
\end{equation}
where $\tau_{p}$ are the interpolation points, and $w_{p}$ are the barycentric weights. By using Lagrange polynomials in the barycentric form, the interpolation weights can be pre-calculated and reused, resulting in computational savings. The interpolation points $\tau_{p}$ are chosen so as to ensure numerical stability, e.g., Chebyshev or Legendre distributions \citep{berrut_barycentric_2004}. We denote the sets of interpolation points $\mathcal{P} := (\mathcal{P}_x, \mathcal{P}_u)$, where  $\mathcal{P}_x := \{\tau_0, \tau_1, ..., \tau_{n_{\mathcal{P}_x}} \}$ and $\mathcal{P}_u := \{\tau_0, \tau_1, ..., \tau_{n_{\mathcal{P}_u}} \}$.

\begin{figure}
    \centering
    \includegraphics[width=0.9\columnwidth]{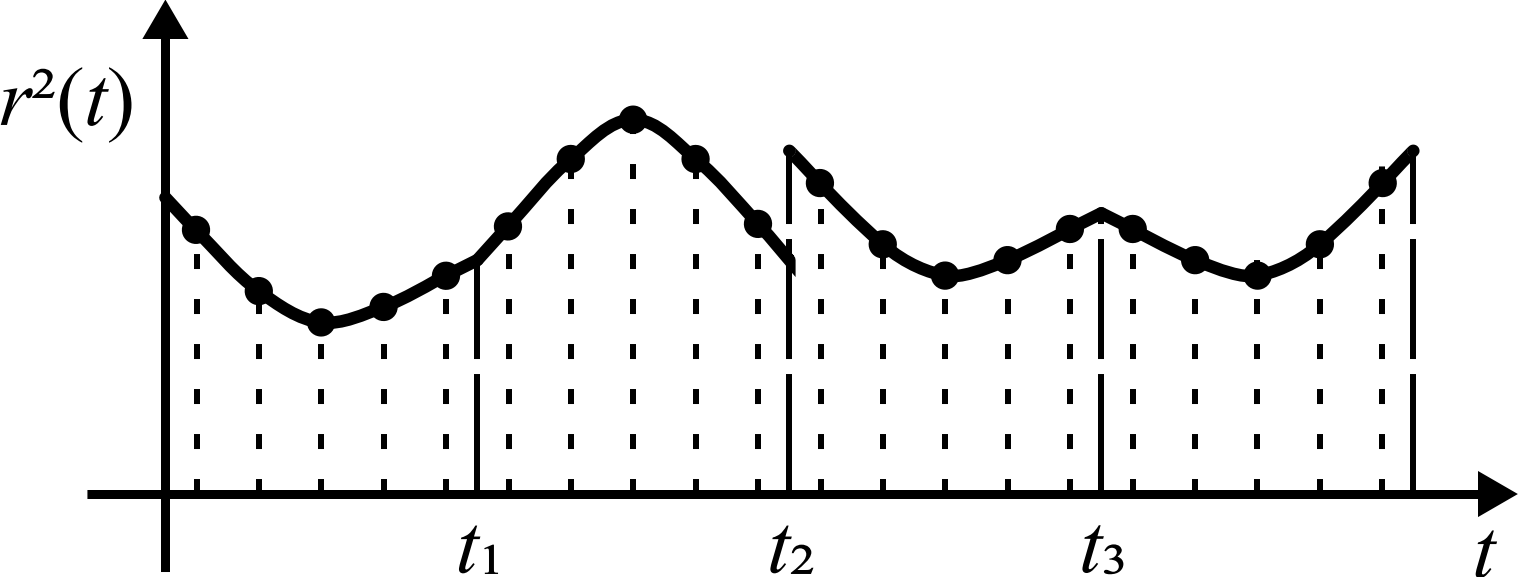}
    \caption{Illustration of Gaussian quadrature with $n_{\mathcal{Q}}=5$.}
    \label{fig:quadrature}
\end{figure}

Thirdly, the mean integrated residual is approximated by a numerical integration scheme. Any quadrature rule can be used, for the case of Gaussian quadrature rule of degree $n_{\mathcal{Q}}$ with quadrature nodes $\tau_q \in [-1, 1]$, the cost functional \eqref{eq:irm} at each interval $i$ is discretised as
\begin{equation}
   \frac{h_i}{2 n_r} \sum_{q=1}^{n_{\mathcal{Q}}} w_{q}
    \bigg{\lVert} r \bigg(
    \frac{2 \dot{\tilde{\boldsymbol{x}}}_i (\tau_{q}, x)}{h_i(t_f - t_0)},
    \tilde{\boldsymbol{x}}_i (\tau_{q}, x),
    \tilde{\boldsymbol{u}}_i (\tau_{q}, u) \bigg) \bigg\rVert ^2,
\label{eq:fn}
\end{equation}
where $\tau_{q}$ are the quadrature points, and $w_{q}$ are the quadrature weights, which can also be pre-computed and reused. Quadrature rules with high accuracy include Gauss-Legendre and Clenshaw-Curtis \citep{trefethen_is_2008}. Figure \ref{fig:quadrature} illustrates the quadrature sampling of the integrand at quadrature points in each sub-interval.
We denote the set of quadrature points $\mathcal{Q} := \{\tau_0, \tau_1, ..., \tau_{n_{\mathcal{Q}}} \}$.

Due to change the of variables from $t$ to $\tau_i$, $\dot{\tilde{\boldsymbol{x}}}$ is appropriately scaled in \eqref{eq:fn}. We can use the sum of \eqref{eq:fn} for each interval as the discretized objective function for the transcribed optimization problems.

\subsection{Sequence of Optimization Problems}
In practice, instead of solving a single high-fidelity approximation of \eqref{eq:p}, we solve a sequence of approximating problems of increasing accuracy. We define the \emph{family} of finite-dimensional optimization problems, parameterized by the mesh $\mathcal{M} := (\mathcal{H}, \mathcal{P}, \mathcal{Q})$, 
\begin{equation}
    \min_{z \in \mathcal{Z}_{\mathcal{M}}} \quad f_{\mathcal{M}}(z),
\tag{$\textbf{P}_{\mathcal{M}}$}
\label{eq:pm}    
\end{equation}
where the pair $z := (x, u)$, $f_{\mathcal{M}} = \sum_{i=1}^{n_{\mathcal{H}}}$ \eqref{eq:fn}, and the family of feasible sets
\begin{equation}
\begin{aligned}
    \mathcal{Z}_{\mathcal{M}} := \{ (x, u) \, | \,
    &x^i_p \in \mathcal{X} \quad \forall p \in \{0, 1, ..., n_{\mathcal{P}_x}\}, \\
    &u^i_p \in \mathcal{U} \quad \forall p \in \{0, 1, ..., n_{\mathcal{P}_u}\}, \\
    &\forall i \in \{1, 2, ..., n_{\mathcal{H}}\},\\
    &x^{i-1}_{n_{\mathcal{P}_x}} = x^i_0 \quad \forall i \in \{2, 3, ..., n_{\mathcal{H}}\} \\
    &x^1_0 \in \mathcal{X}_0, \, x^{n_{\mathcal{H}}}_{n_{\mathcal{P}_x}} \in \mathcal{X}_f \}.
\end{aligned}
\end{equation}
In practice, we wish to \emph{warm-start} each problem with the solution of the previous. Thus we require the constraint sets to be nested, i.e., $\mathcal{Z}_{\mathcal{M}_m} \subseteq \mathcal{Z}_{\mathcal{M}_{m+1}}\, \forall m \in \{1, 2, ...\}$, where each $m$ is a mesh refinement iteration.

The structure of the constraints in \eqref{eq:pm} can be simplified (via elimination of variables for the continuity constraints) to only upper and lower bounds on~$z$. Box-constrained optimization algorithms are tailored for this problem structure. Moreover, algorithms for least-squares objectives may further exploit the problem structure.

A family of optimality functions $\theta_{\mathcal{M}}(\cdot)$ corresponding to problems \eqref{eq:pm} can be used as measures of optimality. An optimality function is a nonpositive scalar-valued function based on optimality conditions of an optimization problem. Minimizers are zeros of a problem's optimality functions \cite[Sect. 1.1.1]{polak_optimization_1997}. For our family of box-constrained problems \eqref{eq:pm}, first-order quadratic approximations can be used as optimality functions, i.e.,
\begin{equation}
    \theta_{\mathcal{M}} (z) := \min_{z' \in \mathcal{Z}_{\mathcal{M}}} \langle \nabla f_{\mathcal{M}}(z), z' - z \rangle  + \frac{1}{2} \|z' - z\|^2,
\end{equation}
which can easily be evaluated if $\nabla f_{\mathcal{M}}(\cdot)$ is known. This measure of optimality can be used to quantify the progress of an iterative algorithm at each step. 

\section{Mesh Refinement}
\label{sec:refinement}
Two strategies to approximate a solution to \ref{eq:p} are described. Both approaches contain an outer loop for mesh refinement, which yields a sequence of transcribed problems $\{ \textbf{P}_{\mathcal{M}_0}, \textbf{P}_{\mathcal{M}_1}, \textbf{P}_{\mathcal{M}_2}, ... \}$. Moreover, an inner loop iterates on the transcribed problems using an optimization algorithm tailored for box constraints. 

Firstly, we propose a \emph{progressive} refinement strategy that uses an early termination test based on the convergence rate of the transcribed problem. A quadrature error test is included to ensure accurate integration. Secondly, a \emph{predictive} refinement strategy is also described, for comparison. Because no mesh refinement scheme has yet been developed for integrated residual methods, we synthesize an $h$-method inspired by the state-of-the-art \citep[Sect. 4.11]{betts_practical_2020}, extended to include the quadrature error test.

\subsection{Refinement Strategies}
Starting with a coarse mesh $\mathcal{M}_0$ and an initial guess $z_0 \in \mathcal{Z}_{\mathcal{M}_0}$, Algorithm~\ref{alg:et} progressively refines the mesh, terminating each transcribed problem when slow convergence is detected (as proposed in \citet[Theorem 3.3.24]{polak_optimization_1997}, under convexity assumptions).
\begin{algorithm}[b]
\caption{Refinement with Early Termination}
\begin{algorithmic}
    \Require $\mathcal{M}_0$, $z_0 \in \mathcal{Z}_{\mathcal{M}_0}$, $\epsilon_{\rho}$, $\epsilon_f$
    \State $\mathcal{M} \gets \mathcal{M}_0$
    \State $k \gets 0$
    \Repeat
        \State $z_{k+1} \gets A_{\ref{alg:bci}}(\mathcal{M}, z_k)$
        \Comment{Use Algorithm \ref{alg:bci}}
        \State $\rho_k \gets \dfrac{f_{\mathcal{M}}(z_{k+1})}{f_{\mathcal{M}}(z_k)}$
        \Comment{Estimated convergence rate}
        \If{$\rho_k \geq \epsilon_{\rho}$}
            \State $k \gets k+1$
        \Else
            \State $(\mathcal{M}, z_k) \gets A_{\ref{alg:progressive}}(\mathcal{M}, z_k)$
            \Comment{Use Algorithm \ref{alg:progressive}}
        \EndIf
    \Until{$f_{\mathcal{M}}(z_k) \leq \epsilon_f$}
\end{algorithmic}
\label{alg:et}
\end{algorithm}
The estimated convergence rate $\rho_k$  is the ratio of integrated residuals between inner loop iterations. Once $\rho_k$ exceeds a convergence threshold $\epsilon_{\rho}$, the optimization is terminated, and Algorithm \ref{alg:progressive} ($A_{\ref{alg:progressive}}$) refines the mesh and interpolates a new $z$ approximation.
\begin{algorithm}[b]
\caption{Progressive Mesh Refinement}
\begin{algorithmic}
    \Require $\mathcal{M}$, $z \in \mathcal{Z}_{\mathcal{M}}$, $\epsilon_q$, $\epsilon_f$
    \State $(\mathcal{H}, \mathcal{P}, \mathcal{Q}) \gets \mathcal{M}$
    \State $e_q \gets \dfrac{|f_{\mathcal{M}}(z) - f_{{\mathcal{M}^+}}(z)|}{1 + f_{{\mathcal{M}^+}}(z)}$
    \Comment{Estimate quadrature error}
    \If{$e_q > \epsilon_q$}
        \State increase $n_{\mathcal{Q}}$
    \Else
        \For{$i \in \{1, 2, ..., n_{\mathcal{H}}\}$}
            \If{$f^i_{\mathcal{M}}(z)/h_i \geq \epsilon_f$}
               \State Partition $h_i \in \mathcal{H}$ into two intervals
            \EndIf
        \EndFor
    \EndIf
    \State $\mathcal{M} \gets (\mathcal{H}, \mathcal{P}, \mathcal{Q})$
    \State Interpolate $z$ into $\mathcal{Z}_{\mathcal{M}}$\\
    \Return $(\mathcal{M}, z)$
\end{algorithmic}
\label{alg:progressive}
\end{algorithm}
An estimate of the quadrature error can be obtained by using a nested quadrature rule. To aid notation, $\mathcal{M}^+$ denotes a mesh with a higher degree quadrature than $\mathcal{M}$. Should the local relative quadrature error $e_q$ exceed a tolerance $\epsilon_q$, the quadrature degree is increased (e.g., $n_{\mathcal{Q}} \gets n_{\mathcal{Q}} + 1$).

The second mesh refinement approach refines the mesh without early termination. \citet{betts_practical_2020} models the differential equation error as being proportional to $h^{p+1}$, where $p$ is the polynomial degree. The synthesized $h$-method predicts the number of partitions as $(e_i/\epsilon_f)^{\frac{1}{p+1}}$, where $e_i := f^i_{\mathcal{M}}(z)/h_i$ and $\epsilon_f$ is the tolerance. The optimization routine is terminated once the optimality function $\theta_{\mathcal{M}}(\cdot)$ reaches an optimality threshold $\epsilon_{\theta}$, as described in Algorithm~\ref{alg:noet}.
\begin{algorithm}[t]
\caption{Refinement without Early Termination}
\begin{algorithmic}
    \Require $\mathcal{M}_0$, $z_0 \in \mathcal{Z}_{\mathcal{M}_0}$, $\epsilon_{\theta}$, $\epsilon_f$
    \State $\mathcal{M} \gets \mathcal{M}_0$
    \State $k \gets 0$
    \Repeat
        \State $z_{k+1} \gets A_{\ref{alg:bci}}(\mathcal{M}, z_k)$
        \Comment{Use Algorithm \ref{alg:bci}}        
        \If{$\theta_{\mathcal{M}}(z) \leq \epsilon_{\theta}$}
            \State $k \gets k+1$
        \Else
            \State $(\mathcal{M}, z_k) \gets A_{\ref{alg:predictive}}(\mathcal{M}, z_k)$
            \Comment{Use Algorithm \ref{alg:predictive}}
        \EndIf
    \Until{$f_{\mathcal{M}}(z_k) \leq \epsilon_f$}
\end{algorithmic}
\label{alg:noet}
\end{algorithm}
Once each transcribed problem is solved, Algorithm~\ref{alg:predictive} ($A_{\ref{alg:predictive}}$) refines the mesh in a predictive manner. While the quadrature branch is similar to Algorithm~\ref{alg:progressive}, the mesh intervals are partitioned into a potentially larger number of intervals. The number of partitions is estimated based on the approximation's degree of polynomial $p$ (e.g., degree of $\mathcal{P}_x$), limited to a maximum of $p_{max}$ partitions. 

\begin{algorithm}[t]
\caption{Predictive Mesh Refinement}
\begin{algorithmic}
    \Require $\mathcal{M}$, $z \in \mathcal{Z}_{\mathcal{M}}$, $\epsilon_q$, $\epsilon_f$, $p$, $p_{max}$
    \State $(\mathcal{H}, \mathcal{P}, \mathcal{Q}) \gets \mathcal{M}$
    \State $e_q \gets \dfrac{|f_{\mathcal{M}}(z) - f_{{\mathcal{M}^+}}(z)|}{1 + f_{{\mathcal{M}^+}}(z)}$
    \Comment{Estimate quadrature error}
    \If{$e_q > \epsilon_q$}
        \State increase $n_{\mathcal{Q}}$
    \Else
        \For{$i \in \{1, 2, ..., n_{\mathcal{H}}\}$}
            \State $e_i \gets f^i_{\mathcal{M}}(z)/h_i$
            \If{$e_i \geq \epsilon_f$}
                \State Partition $h_i \in \mathcal{H}$ into
                \State $\min\{ \lceil (e_i/\epsilon_f)^{\frac{1}{p+1}} \rceil, p_{max} \}$ intervals
            \EndIf
        \EndFor
    \EndIf
    \State $\mathcal{M} \gets (\mathcal{H}, \mathcal{P}, \mathcal{Q})$
    \State Interpolate $z$ into $\mathcal{Z}_{\mathcal{M}}$\\
    \Return $(\mathcal{M}, z)$
\end{algorithmic}
\label{alg:predictive}
\end{algorithm}

\subsection{Optimization Algorithm}
Many optimization methods are compatible with Algorithms~\ref{alg:et} and \ref{alg:noet}. Because the family of optimization problems \ref{eq:pm} can be simplified into only having box constraints, projected gradient methods \citep{bertsekas_projected_1982} can be used. This class of methods is favorable because the constraints are enforced at all iterations, can easily be warm-started, and are easy to implement.

For a set of simple bounds $\mathcal{Z} = \{ z_{\ell} \leq z \leq z_u \}$, the projection operator $P_{\mathcal{Z}}(\cdot)$ is defined as
\begin{equation}
    P_{\mathcal{Z}}(z^j) := 
    \begin{cases}
        z_u^j \quad & \textrm{if} \,z^j > z_u^j, \\
        z_{\ell}^j  & \textrm{if} \,z^j < z_{\ell}^j, \\
        z^j         & \textrm{otherwise},
    \end{cases}
\end{equation}
for each $j^{\text{th}}$ index of $z$. The set of binding indices $\mathcal{B}$ corresponds to active (or almost active) bound constraints, defined as
\begin{multline}
    \mathcal{B}(z) := \big\{ j \,|\, z^j_{\ell} \leq z^j \leq z^j_{\ell} + \epsilon, \, \nabla f(z)^j > 0 \\
    \textrm{or} \quad z^j_u - \epsilon \leq z^j \leq z^j_u, \, \nabla f(z)^j < 0 \big\},
\end{multline}
where $\epsilon = \textrm{min}(\epsilon_z, \| z - P_{\mathcal{Z}}(z - \nabla f(z)) \|)$.

Algorithm \ref{alg:bci} describes an iteration of a projected gradient method.
\begin{algorithm}[tb]
\caption{Box-constrained iteration}
\begin{algorithmic}
\Require $\mathcal{M}$, $z_0$
\State $m \gets 0$
\Repeat
    \State $d_m$, $\alpha_m$ from \citet{bertsekas_projected_1982}
    \State $z_{m+1} \gets P_{\mathcal{Z}_{\mathcal{M}}}(z_m + \alpha_m d_m)$
    \State $m \gets m+1$
\Until{$\mathcal{B}(z_{m+1}) = \mathcal{B}(z_m)$}\\
\Return $z_{m+1}$
\end{algorithmic}
\label{alg:bci}
\end{algorithm}
\citet{bertsekas_projected_1982} describes how a descent direction $d$ can be obtained from steepest-descent, conjugate-gradient, and Newton directions. An Armijo-like line search method is also proposed to compute the step size~$\alpha$.

\section{Example}
\label{sec:example}
We assess the efficiency of the mesh refinement strategies by applying them to a nonlinear problem, whose solution cannot be captured exactly by piece-wise polynomials.

In this problem \cite[Sect. 6]{kelly_introduction_2017}, a simple pendulum is attached to a cart. The dynamics are modelled by the differential equations
\begin{multline}
    \ddot{q}_1      = \frac{u_1 + m_2 \ell \sin(q_2) \dot{q}_2^2 + m_2 g \cos(q_2) \sin(q_2)}{m_1 + m_2 \sin^2(q_2)}, \, \ddot{q}_2 = \\
    \frac{u_1 \cos(q_2) +m_2 \ell \cos(q_2) \sin(q_2) q_2^2 + (m_1 + m_2) g \sin(q_2)}{m_1 + m_2 \sin^2(q_2)},
\label{eq:cartpole}
\end{multline}
where $q_1$ and $\dot{q}_1$ are the cart's position (m) and velocity (m/s), $q_2$ and $\dot{q}_2$ are the pendulum's angle (rad) and angular rate (rad/s), and $u_1$ is the control input (N). The masses for the cart and pendulum are set to $m_1 = 1$\,kg and $m_2 = 0.3$\,kg, the pendulum's length $\ell = 0.5$\,m, and gravity's acceleration $g = 9.81$\,m/s$^2$. Starting from rest
\begin{equation}
\begin{aligned}
    q_1(0) &= 0 \quad \quad q_2(0) = 0\\
    \dot{q}_1(0) &= 0 \quad \quad \dot{q}_2(0) = 0,
\end{aligned}
\label{eq:x0}
\end{equation}
the cart must be actuated so as to swing up the pendulum in two seconds, i.e., satisfying
\begin{equation}
\begin{aligned}
    q_1(2) &= 1 \quad \quad q_2(2) = \pi\\
    \dot{q}_1(2) &= 0 \quad \quad \dot{q}_2(2) = 0.
\end{aligned}
\label{eq:xf}
\end{equation}

The swing-up example is formulated as \eqref{eq:dfp} with $q_1$, $\dot{q}_1$, $q_2$, $\dot{q}_2$ as differential dynamic variables, and $u_1$ as an algebraic dynamic variable. The differential equations~\eqref{eq:cartpole} are included as residuals in $r(\cdot, \cdot, \cdot)$, and boundary conditions~\eqref{eq:x0}, \eqref{eq:xf} are enforced as initial and terminal constraints.

The DFP is discretised with piece-wise quadratic polynomials ($n_{\mathcal{P}_x}=3$, $n_{\mathcal{P}_u}=3$), and is solved using a projected conjugate gradient calling an Armijo-like line-search \citep{schwartz_family_1997}.
We solve the DFP until an integrated residual tolerance $\epsilon_f = 1.15 \times 10^{-1}$ is reached. We compare three approaches:
\begin{enumerate}
    \item No mesh refinement, using a sufficiently large number of intervals $n_{\mathcal{H}} = 20$, and quadrature $n_{\mathcal{Q}} = 4$;
    \item The refinement scheme without early termination (Algorithm \ref{alg:noet}), starting with $n_{\mathcal{H}} = 2$, $n_{\mathcal{Q}} = 4$, $\epsilon_{q} = 1\times10^{-1}$, $\epsilon_{\theta} = 2\times10^{-4}$, and $p_{max} = 4$;
    \item The refinement scheme with early termination (Algorithm \ref{alg:et}) starting with $n_{\mathcal{H}} = 2$, $n_{\mathcal{Q}} = 4$, $\epsilon_{q} = 1\times10^{-1}$, and $\epsilon_{\rho} = 0.999$.
\end{enumerate}

To ensure a fair comparison, the values of $n_{\mathcal{H}}$, $\epsilon_{\theta}$, and $\epsilon_{\rho}$ were experimentally chosen to yield the least number of $\nabla r$ evaluations for each approach.

The three strategies were implemented in the Julia programming language, solved using an Intel\textregistered\, Core\texttrademark\, i7-11370H laptop running Microsoft Windows 10.0.22621 with 16 GB of RAM. A summary of the results is presented in Table \ref{tab:results}, where solve times are the minimum of 10 trials.
\begin{table}
\caption{Summary of performance results}
\label{tab:results}
\begin{center}
\begin{tabular}{|c||c|c|}
\hline
                            & Total $\nabla r$ evaluations & Solve time\\
\hline
No mesh refinement          & 44 960                & 4.97 s\\
\hline
MR, no early termination    & 40 900               & 4.77 s\\
\hline
MR with early termination   & 10 088                & 1.23 s\\
\hline
\end{tabular}
\end{center}
\end{table}

\begin{figure}
    \centering
    \includegraphics[width=\columnwidth]{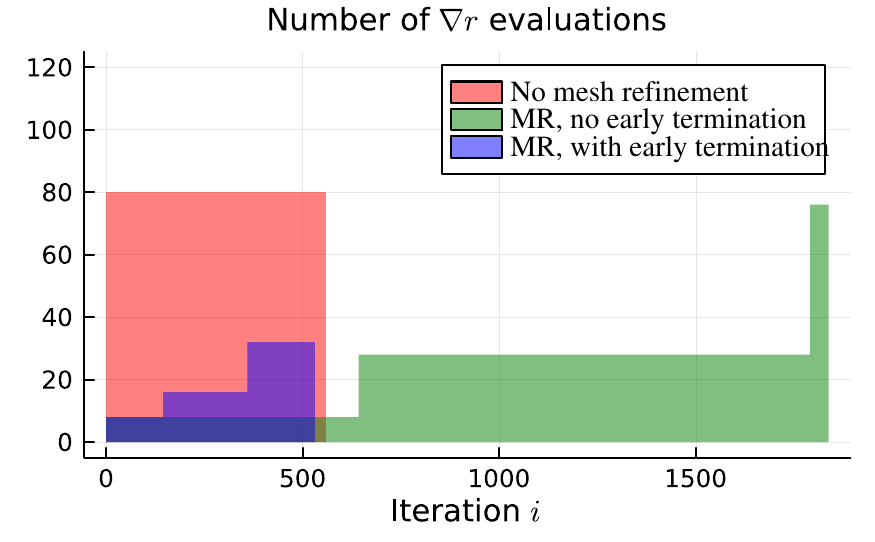}
    \caption{Performance comparison based on $\nabla r$ evaluations}
    \label{fig:calls}
\end{figure}

Figure \ref{fig:calls} shows the number of evaluations of the residuals Jacobian function $\nabla r$ at each iteration of the optimization algorithm.
Evaluations of $\nabla r$ dominate the total computational burden, with mesh refinement and interpolation procedures becoming insignificant. The colored area in each plot is therefore representative of the overall performance.

\begin{figure}
    \centering
    \includegraphics[width=\columnwidth]{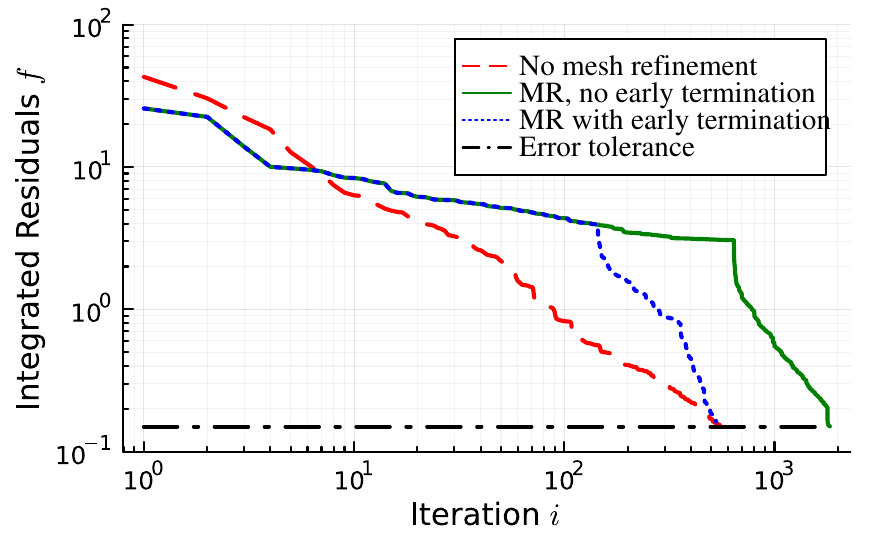}
    \caption{Convergence of the objective function $f$ for each strategy. Mesh refinements appear as sharp kinks.}
    \label{fig:residuals}
\end{figure}

Figure \ref{fig:residuals} shows the mean integrated residual $f$ against the iteration number on a log-log plot. 
The convergence of different approaches can be examined.

The performance of the mesh refinement strategies can be described by a joint analysis of Figures \ref{fig:calls} and \ref{fig:residuals}. Starting with a coarse mesh is advantageous not only in the computation per iteration (smaller optimization problem), but also allows for faster convergence on the first iterations.

Each optimization problem need not be solved to high precision, especially initially. A phenomenon of diminishing returns is observed, where iterating a small optimization problem no longer provides a significant decrease in $f_{\mathcal{M}}$, compared to iterating on a larger problem, thus justifying early termination.

\section{Conclusions}
\label{sec:conclusions}
We have described a novel mesh refinement strategy for solving DFPs using the general class of integrated residual methods. On the example problem, the progressive refinement was shown to outperform the conventional predictive approach by up to a factor of three. Early termination of the optimization problems was shown to improve the efficiency of the overall numerical method.

The following are suggestions for future research directions:
\begin{itemize}
    \item A more sophisticated mesh refinement procedure: an $h$-method is modest compared to state-of-the-art ($hp$ or $ph$) methods. Care must be taken on how the choice between $h$ and $p$ refinements is made, particularly considering the sparsity structure of the transcribed optimization problem;
    \item A selection of matured optimization algorithms: Both first-order and second-order methods should be explored. Projected search methods may be paired with recently proposed line-search algorithms that provide stronger, Wolfe-like conditions \citep{ferry_projected-search_2021};
    \item Optimized \emph{hyper-parameters} for Algorithms \ref{alg:et}, \ref{alg:progressive}, and \ref{alg:bci}: Zero-order methods such as Bayesian optimization may be implemented to optimize the performance and/or robustness of the algorithms;
    \item A comparison between integrated residuals methods, and state-of-the-art direct collocation and multiple shooting methods;
    \item An extension of the refinement method to dynamic optimization problems, possibly as per \citet{CDCInteresso}.
\end{itemize}


\bibliography{references}             
                                                   







\end{document}